\documentclass[12pt,a4paper]{amsart}
\usepackage{amsmath,amssymb}
\usepackage{graphicx}
\usepackage{amsfonts}
\newtheorem{Lma}{Lemma}
\newtheorem{Thm}{Theorem}
\newtheorem{Prop}{Proposition}

\theoremstyle{remark}
\newtheorem{Rmk}{Remark}

\begin{document}
\title[Volume Estimates of Expanding Gradient Ricci Solitons]
{Volume Estimates And The Asymptotic Behavior of Expanding Gradient Ricci Solitons}
\author{Chih-Wei Chen}
\address{Department of Mathematics, National Taiwan University, Taiwan}
\address{Institut Fourier, Universit\'e Joseph Fourier, France}
\email{BabbageTW@gmail.com}
%\thanks{Thanks for Author One.}
\date{April 21, 2011}
%\keywords{Ricci soliton, injectivity radius, curvature decay}

\begin{abstract}
We study the asymptotic volume ratio of non-steady gradient Ricci solitons.
Moreover, a local estimate of the volume ratio is obtained for
expanding solitons which satisfy $\lim_{dist(O,x)\rightarrow\infty} |Sect|\cdot dist(O,x)^2=0$.
Therefore, for such a soliton, we can show that it must have $\mathbb{R}^n$ as one of its tangent cone at infinity.
(Here we assume that the soliton is simply connected at infinity, has only one end and $n\geq 3$.)
\end{abstract}

\maketitle

%\subjclass[2000]{Primary 05C38, 15A15; Secondary 05A15, 15A18}

%%%%%%%%%%%%%%%%%%%%%%%%%%%%%%%%%%%%%%%%%%%%%%%%%%%%%%%%%%%%%%%%%%%%%%%%%%%%%%%%%%%%%%%%%%%%%%%%%%%%%
\section{Introduction}%%%%%%%%%%%%%%%%%%%%%%%%%%%%%%%%%%%%%%%%%%%%%%%%%%%%%%%%%%%%%%%%%%%%%%%%%%%%%%%
%%%%%%%%%%%%%%%%%%%%%%%%%%%%%%%%%%%%%%%%%%%%%%%%%%%%%%%%%%%%%%%%%%%%%%%%%%%%%%%%%%%%%%%%%%%%%%%%%%%%%

The Ricci solitons, which are generalizations of the Einstein manifolds, are important solutions
to the Ricci flow. Besides the advantage of having explicit equations, they occur in the analysis
of blow-up limits near singularities. In this article, we only discuss the complete non-compact
solitons, which are much more complicated than the compact ones.
In the three-dimensional case, the classification of shrinking solitons under some reasonable
conditions leads to the performance of surgery. For higher dimensional cases, some results about the
classification of solitons were obtained in the last four years, e.g. \cite{ Chen, Na, NW, PW, CM, Zh}.
These results were derived under various curvature assumptions such as locally conformally flat,
constant scalar curvature, nonnegative Ricci curvature(for expanding solitons) or bounded nonnegative
curvature operator(for shrinking solitons when $n=4$.) In this article, we try to understand the geometry
of solitons which are not Ricci-nonnegative.

Besides the studies on the classification, there are some results and conjectures about the non-existence.
We recall some classical non-existence theorems about general complete Riemannian manifolds.
In \cite{BKN}, S. Bando, A. Kasue and H. Nakajima proved that there exists no manifold with $Ric\geq 0$,
$\lim_{dist(O,x)\rightarrow\infty} |Sect|\cdot dist(O,x)^2=0$ and $Vol(B_s)\geq Cs^n$ for all geodesic balls $B_s$ with
radius $s$ and center $O$, where we use $C$ to denote various constants.
Another non-existence result due to R. E. Greene and H. Wu \cite{GW} and G. Drees \cite{D} states that
there exists no manifold with positive sectional curvature and $\lim_{dist(O,x)\rightarrow\infty} R\cdot dist(O,x)^2=0$
except for $n=4$ or $8$. Here, and afterwards, $R$ always stands for the scalar curvature.
An approach to achieve these non-existence results is to study the tangent cones at infinity of such manifolds.
Indeed, we prove that if a non-flat non-steady Ricci soliton $M$ satisfies
$\lim_{dist(O,x)\rightarrow\infty} |Sect|\cdot dist(O,x)^2=0$ and is simply connected at infinity,
then each tangent cone at infinity of $M$ is the Euclidean space $\mathbb{R}^n$.
(Here we assume that the soliton has only one end and has dimension $n\geq 3.$)
For the case of gradient expanding solitons, B.-L. Chen and X.-P. Zhu \cite{ChZh} proved that
they cannot have $Sect\geq 0$ and positively $\epsilon$-pinched Ricci curvature, i.e. $Ric\geq \epsilon Rg$ and $R>0$, when $n\geq 3$.
It is still unknown that whether the condition on the sectional curvature can be discarded or not.
(See also \cite{MC,M} and \cite{N}, section 3, for some discussions about this.)

We should mention that there were other results about the non-existence of Ricci solitons. For example,
S. Pigola, M. Rimoldi and A. G. Setti \cite{PRS} proved that only trivial expanding solitons
$(M,g,f\equiv const.)$ can satisfy $|\nabla f|\in L^p(M,e^{-f}dvol)$ for some $1\leq p\leq\infty$.
They also proved that an expanding soliton must be flat provided that $0\leq R\in L^1(M,e^{-f}dvol)$.
Note that, up to today, the growth of $f$ is unknown for expanding solitons unless we have some control
on the Ricci curvature.
On the other hand, A. Deruelle \cite{De} proved the rigidity of steady Ricci solitons with
$Sect\geq 0$ and $R\in L^1(M)$.

In the next section, we study the behavior of $f$ on expanding solitons with $|Ric|\leq C\cdot dist(O,x)^{-\varepsilon}.$
In Section 3, we derive various bounds of the asymptotic volume ratio for non-steady solitons with different
decay types of curvature. A uniform volume bound for all geodesic balls are derived under the condition that
$\lim_{dist(O,x)\rightarrow\infty} |Sect|\cdot dist(O,x)^2=0$.
Use this estimate, we can study the tangent cones at infinity of non-steady Ricci solitons in the last section.
Throughout this article,
$s$ denotes either the distance $s(x):=dist(O,x)$ or the radius of a geodesic ball with center $O$
and $C$ denotes some positive constant which may vary in different contexts.

The author appreciates the hospitality and suggestions of professors Xi-Ping Zhu and Bing-Long Chen
when he visited Guangzhou and also the discussions with Hui-Ling Gu and Zhu-Hong Zhang.
He would like to thank his advisors, G\'erard Besson and Yng-Ing Lee, for their encouragements and
discussions.
Special thanks to professor Gilles Carron for his enthusiasm on this article.
Surely the author cannot finish this article without the help of Carron.

Some inaccurate quotes about \cite{PRS} has been corrected in this version, thanks to the communication from Michele Rimoldi.

%%%%%%%%%%%%%%%%%%%%%%%%%%%%%%%%%%%%%%%%%%%%%%%%%%%%%%%%%%%%%%%%%%%%%%%%%%%%%%%%%%%%%%%%%%%%%
\section{The growth of $f$ on Ricci solitons with certain curvature decay}
%%%%%%%%%%%%%%%%%%%%%%%%%%%%%%%%%%%%%%%%%%%%%%%%%%%%%%%%%%%%%%%%%%%%%%%%%%%%%%%%%%%%%%%%%%%%%

Let $(M,g,f)$ be a complete non-compact expanding gradient Ricci soliton, which satisfies the following equation:
$$R_{ij}+\nabla_i\nabla_jf=-g_{ij},$$ and $R$ be the scalar curvature of $(M,g,f)$. The following
three lemmas are well-known.
{\Lma {\rm (R. S. Hamilton, \cite{Ha})} We have $R+|\nabla f|^2+2f=C_1$ for some constant $C_1$
which can be absorbed by $f$.}
{\Lma {\rm (R. S. Hamilton, \cite{Ha93})} The time-independent Harnack quantity
$\Delta R -\left<\nabla R, \nabla f\right>+2(R+|Ric|^2)$
vanishs on $(M,g,f)$.}
{\Lma {\rm (B.-L. Chen, \cite{Chen})} We have $R\geq-C_2$ for some constant $C_2>0$.}\\

As we can see in Lemma 1, there is a normalization on the function $f$ which is usually
used to simplify proofs in many cases.
Considering the derivative of $f$, we have the following a priori relation between
$\nabla f$ and $\nabla R$.

{\Thm%%%%%%%%%%%%%%%%%%%%%%%%%%%%%%%%%%%%%%%%%%%%%%%%%%%%%%%%%%%%%%%%%%%%%%
If $\nabla f(p) =0$ for some $p\in M$, then $\nabla R(p)=0$. On the other hand,
if $\nabla R(p)=0$ and $R(p)<-\frac{n-1+|Ric|^2}{2}$, then $\nabla f(p)=0$.
}%%%%%%%%%%%%%%%%%%%%%%%%%%%%%%%%%%%%%%%%%%%%%%%%%%%%%%%%%%%%%%%%%%%%%%%%%%

\begin{proof}
If $\nabla f(p) =0$, then we have $d R = 2Ric(\nabla f,\cdot)=0$ at $p$.

On the other hand, suppose $\nabla R(p)=0$ and $\nabla f(p)\neq 0$,
we claim that $R(p)\geq -\frac{n-1+|Ric|^2}{2}$. Since $|\nabla f|$ is locally Lipschitz,
by Lemma 1, we have
$$-2\nabla f=\nabla |\nabla f|^2 = 2 |\nabla f|\cdot \nabla |\nabla f|.$$
By using the assumption $\nabla f(p)\neq 0$ and Kato's inequality $|\nabla |\nabla f||\leq |Hess f|$, we can divide both sides by $|\nabla f|$ and get
$$ |g+Ric|^2=|Hess f|^2\geq |\nabla |\nabla f||^2=\left| \frac{\nabla f}{|\nabla f|} \right|^2=1,$$ i.e.
$n+2R+|Ric|^2\geq 1$. Hence $R\geq -\frac{n-1+|Ric|^2}{2}$ at $p$.
The statement of theorem follows by reduction to the absurd.

\end{proof}

{\Rmk Similar computation holds for shrinking solitons.}\\

Given a fixed point $O\in M$, we set $s=dist(O,x)$ and $\gamma(s)$ be a unit-speed
minimizing geodesic connecting $O$ and $x$, where $x\in M$ is chosen arbitrarily.
We use the notation $'$ to denote the differentiation with respect to $s$ along $\gamma(s)$.
The following proposition, which seems to appear first time in the
literature in \cite{Z}, is an easy consequence of Lemmas 1 and 2.
{\Prop%%%%%%%%%%%%%%%%%%%%%%%%%%%%%%%%%%%%%%%%%%%%%%%%%%%%%%%%%%%%%%%%%%
For every expanding soliton $M$, we have $|f'(x)|\leq |\nabla f(x)|\leq s+L(O)$,
where $L(x)=\sqrt{C_1+C_2-2f(x)}=\sqrt{C_2+R(x)+|\nabla f(x)|^2}$.
Moreover, when $Ric\geq 0$, we have $f'(x)\leq -s+f'(O)$.
}%%%%%%%%%%%%%%%%%%%%%%%%%%%%%%%%%%%%%%%%%%%%%%%%%%%%%%%%%%%%%%%%%%%%%%%

\begin{proof}
Since $-C_2+|\nabla f|^2+2f \leq R+|\nabla f|^2+2f=C_1$, we have $|\nabla f| \leq \sqrt{C_1+C_2-2f}=L.$
Combining with $\nabla L=\frac{-\nabla f}{\sqrt{C_1+C_2-2f}}$, we have $|\nabla L|\leq 1.$

Integrating it from the point $O$ to some point $x=\gamma(s)$ along $\gamma$,
we have $$L(x)-L(O)= \int_0^s L'\leq  \int_0^s |\nabla L| \leq s.$$
Hence, $|f'(x)|\leq |\nabla f(x)| \leq L(x) \leq s+L(O).$

When $Ric\geq 0$, $$\int_0^s f'' \leq \int_0^s Ric(\gamma',\gamma') + \int_0^s f'' =-s$$ implies that
$f'(x)\leq -s+f'(O).$
\end{proof}

From this proposition, it is easy to see that for every expanding gradient Ricci soliton which has
nonnegative Ricci curvature, the potential function $f(x)$ must decrease quadratically in $s$.
The following theorem shows that this property holds for expanding solitons whose Ricci curvatures may be negative.

{\Thm%%%%%%%%%%%%%%%%%%%%%%%%%%%%%%%%%%%%%%%%%%%%%%%%%%%%%%%%%%%%%%%%%%%%%%%%%%%%%%%%%%%%%%%%%%%%%%%%%%%%%
If $|Ric|\leq Cs^{-\varepsilon}, s\equiv dist(O,x),$ for some constant $\varepsilon<1$ and some point $O\in M$,
then there exists a point $p\in M$ and $C_3,C_4>0$ such that $|Ric|\leq C_3\cdot dist(p,x)^{-\varepsilon}$
and $f$ satisfies
$$-r\left(1+\frac{C_4}{r^\varepsilon}\right) \leq  f'(x)  \leq -r\left(1-\frac{C_4}{r^\varepsilon}\right),$$
where $r=dist(p,x).$ As a consequence, we have
$$-\frac{1}{2}r^2\left(1+\frac{C_5}{r^\varepsilon}\right) +f(p)
\leq  f(x) \leq  -\frac{1}{2}r^2\left(1-\frac{C_5}{r^\varepsilon}\right) +f(p).$$
}%%%%%%%%%%%%%%%%%%%%%%%%%%%%%%%%%%%%%%%%%%%%%%%%%%%%%%%%%%%%%%%%%%%%%%%%%%%%%%%%%%%%%%%%%%%%%%%%%%%%%%%%%%%

\begin{proof}
From
$$    -C\int_0^s s^{-\varepsilon} + \int_0^s f''
  \leq \int_0^s Ric(\gamma',\gamma')  + \int_0^s f''  =   \int_0^s-1
  \leq C\int_0^s s^{-\varepsilon} + \int_0^s f'' ,$$
we have
$$ -s-C\int_0^s s^{-\varepsilon}\leq \int_0^s f'' \leq  -s+C\int_0^s s^{-\varepsilon}$$
and hence
$$-s\left(1+\frac{C_4}{s^\varepsilon}\right)+f'(O)  \leq  f'(x)  \leq -s\left(1-\frac{C_4}{s^\varepsilon}\right)+f'(O).$$

In order to achieve the conclusion, it is enough to show that $f$ has a critical point $p$
(and then repeat the calculation above.)
This can be observed by considering the geodesic sphere $\partial B_s(O)$ with $s$ very large.
Since $\nabla f\cdot\nabla s$ is negative on such sphere, $\nabla f$ must point inwards.
So $\nabla f=0$ at some point $p$ inside the ball $B_s(O).$

\end{proof}

%  Throughout this article, we use $r$ to denote the distance counting from a critical point $p$ of $f$
%  and $s$ to denote the distance counting from an arbitrary point $O\in M$. It only makes difference in
%  the proof of Corollary 1, which uses the theorem above.

We recall that the potential function grows quadratically on every shrinking gradient Ricci soliton.
This was proved by H.-D. Cao and D.T. Zhou in \cite{CZ}.
Moreover, by using the same proof in Theorem 1,
we have
$r\left(1-\frac{C_4}{r^\varepsilon}\right) \leq  f'(x)  \leq r\left(1+\frac{C_4}{r^\varepsilon}\right)$
and
$\frac{1}{2}r^2\left(1-\frac{C_5}{r^\varepsilon}\right) +f(p)
\leq  f(x)  \leq\frac{1}{2}r^2\left(1+\frac{C_5}{r^\varepsilon}\right) +f(p)$ for shrinking solitons
which satisfy $R_{ij}+\nabla_i\nabla_jf=g_{ij}$ and $|Ric|\leq C\cdot dist(O,x)^{-\varepsilon}$.\\

{\Rmk The condition $|Ric|\leq C s^{-\varepsilon}$ in Theorem 1 can be replaced by
$|Ric(\gamma',\gamma')|\leq C s^{-\varepsilon}$ for all $\gamma$ starting from $O$.
It is worthy to distinguish these two conditions because a cigar-like manifold may
satisfy the second condition while breaks the first one.
}
%%%%%%%%%%%%%%%%%%%%%%%%%%%%%%%%%%%%%%%%%%%%%%%%%%%%%%%%%%%%%%%%%%%%%%%%%%%%%%%%%%%%%%%%%%%%%%%%%%%%%%%
\section{Volume estimates of non-steady gradient Ricci solitons}%%%%%%%%%%%%%%%%%%%%%%%
%%%%%%%%%%%%%%%%%%%%%%%%%%%%%%%%%%%%%%%%%%%%%%%%%%%%%%%%%%%%%%%%%%%%%%%%%%%%%%%%%%%%%%%%%%%%%%%%%%%%%%%
It was mentioned in \cite{CLN} that a complete non-compact expanding gradient Ricci
soliton with $Ric>0$ must have positive asymptotic volume ratio, which was
proved by Hamilton. That is, the limit $\lim_{s\rightarrow\infty} \frac{Vol(B_{s})}{s^{n}}$ exists and is positive.
Indeed, one only need to show $\liminf_{s\rightarrow\infty} \frac{Vol(B_{s})}{s^{n}}>0$ because the upper bound comes
from Bishop-Gromov comparison.
However, when the soliton is not Ricci-nonnegative, the limit may not exist.
In \cite{CN}, J. A. Carrillo and L. Ni proved that $\liminf_{s\rightarrow\infty} \frac{Vol(B_{s})}{s^{n}}>0$
by only assuming that the scalar curvature is nonnegative. We now can weaken the
curvature condition to be $\frac{1}{Vol(B_{s})}\int_{B_{s}}R\geq -Cs^{-\varepsilon },$
where $B_s\subset M$ always denotes the geodesic ball with central point $O$ and radius $s$.
We also derive an upper bound estimate for expanding solitons with
$Ric\geq  -Cs^{-\varepsilon}g$ and $\frac{1}{Vol(B_s)}\int_{B_s} R\leq Cs^{-\varepsilon}$
by using the same method.

{\Thm%%%%%%%%%%%%%%%%%%%%%%%%%%%%%%%%%%%%%%%%%%%%%%%%%%%%%%%%%%%%%%%%%%
Let $(M,g,f)$ be a complete non-compact expanding gradient Ricci soliton with scalar curvature $R$.
If there exists $O\in M$ such that $\frac{1}{Vol(B_s)}\int_{B_s} R\geq -Cs^{-\varepsilon}$, where $\varepsilon >0$ is a
constant, then $\liminf_{s\rightarrow\infty} \frac{Vol(B_{s})}{s^{n}} \geq \eta $.
Moreover, if we have $Ric\geq  -Cs^{-\varepsilon}g$ and $\frac{1}{Vol(B_s)}\int_{B_s} R\leq Cs^{-\varepsilon}$,
then $$C^{-1}s^n \leq Vol(B_{s}) \leq Cs^n$$
holds for all $s\geq A,$ where $A$ is a large constant.
}%%%%%%%%%%%%%%%%%%%%%%%%%%%%%%%%%%%%%%%%%%%%%%%%%%%%%%%%%%%%%%%%%%%%%%

\begin{proof}
Taking the trace of the soliton equation $R_{ij}+\nabla _{i}\nabla _{j}f=-g_{ij}$
and integrating it on $B_{s},$ we have%
$$-nVol(B_{s}) =\int_{B_s} R+\int_{B_s} \Delta f
  =  \int_{B_s} R+\int_{\partial B_{s}}\nabla f\cdot \nabla s
\geq \int_{B_s} R-\int_{\partial B_{s}}(s+L(O))$$$$
  =  \int_{B_s} R-(s+L(O))Area(\partial B_{s})
  =  \int_{B_s} R-(s+L(O))\frac{d}{ds}Vol(B_{s}).$$
Therefore,
\begin{eqnarray*}
\frac{d}{ds}\log Vol(B_{s}) &\geq &\frac{1}{(s+L(O))Vol(B_{s})}\int_{B_s} R+\frac{n}{s+L(O)} \\
&=&\frac{1}{(s+L(O))Vol(B_{s})}\int_{B_s} R+\frac{d}{ds}\log (s+L(O))^{n}
\end{eqnarray*}
\begin{eqnarray*}
&\Rightarrow &\frac{d}{ds}\log \frac{Vol(B_{s})}{(s+L(O))^{n}}\geq \frac{1}{%
(s+L(O))Vol(B_{s})}\int_{B_s} R\geq \frac{-C}{(s+L(O))s^{\varepsilon }}\geq \frac{-C}{s^{1+\varepsilon}}\\
&\Rightarrow &\log \frac{Vol(B_{s})}{(s+L(O))^{n}}\geq \int_{\rho }^{s}\frac{-C%
}{s^{1+\varepsilon }}+\log \frac{Vol(B_{\rho })}{(\rho +L(O))^{n}}
=\frac{C}{\varepsilon }s^{-\varepsilon}-\frac{C}{\varepsilon }\rho^{-\varepsilon }
+\log \frac{Vol(B_{\rho })}{(\rho +L(O))^{n}} \\
&& \mbox{ for any positive constant }\rho<s\\
&\Rightarrow &\frac{Vol(B_{s})}{(s+L(O))^{n}}
\geq \left( e^{\frac{C}{\varepsilon }s^{-\varepsilon }-\frac{C}{%
\varepsilon }\rho^{-\varepsilon }}\right) \frac{Vol(B_{\rho })}{(\rho+L(O))^{n}}
\geq e^{-\frac{C}{\varepsilon }\rho^{-\varepsilon }}\cdot \frac{Vol(B_{\rho })}{(\rho+L(O))^{n}}.
\end{eqnarray*}

Hence,
\begin{equation*}
\liminf_{s\rightarrow \infty }\frac{Vol(B_{s})}{s^{n}}\geq
 e^{-\frac{C}{\varepsilon }\rho^{-\varepsilon }}\cdot \frac{Vol(B_{\rho })}{(\rho+L(O))^{n}}\equiv\eta >0.
\end{equation*}

For the reader's convenience, we write down the proof of the upper bound estimate
although it is almost the same to the above one.

From the lower bound of the Ricci curvature and Theorem 2, we have $f'(x)\leq -s+C$ for $s$ large enough.
Together with the lower bound of the averaged scalar curvature, we have
$$-nVol(B_{s}) =\int_{B_s} R+\int_{\partial B_{s}}\nabla f\cdot \nabla s
\leq Cs^{-\varepsilon}Vol(B_s) -(s-C)\frac{d}{ds}Vol(B_s)$$
$$\Rightarrow \frac{-n}{s-C}\leq \frac{C}{(s-C)s^{\varepsilon}} -\frac{d}{ds}\log Vol(B_s)
\leq \frac{C}{s^{1+\frac{\varepsilon}{2}}} -\frac{d}{ds}\log Vol(B_s).$$
Hence we get a similar inequality $\frac{d}{ds}\log \frac{Vol(B_{s})}{(s+C)^{n}}
\leq \frac{C}{s^{1+\frac{\varepsilon}{2}}}.$ The rest of the proof is easy to work out.

\end{proof}

For shrinking gradient Ricci solitons, the same calculation gives the following theorem.
{\Thm%%%%%%%%%%%%%%%%%%%%%%%%%%%%%%%%%%%%%%%%%%%%%%%%%%%%%%%%%%%%%%%%%%
Let $(M,g,f)$ be a complete non-compact shrinking gradient Ricci soliton which satisfies $R_{ij}+\nabla_i\nabla_jf=g_{ij}$.
If there exists $O\in M$ such that $\frac{1}{Vol(B_s)}\int_{B_s} R\leq Cs^{a},$ where $a$ is a nonzero
constant, then its volume ratio $\frac{Vol(B_s)}{s^n}$ is bounded from below by $C\cdot e^{\frac{-1}{a}s^a}$ for $s$ large enough.
When $\frac{1}{Vol(B_s)}\int_{B_s} R\leq \delta_1<n $ (for $s$ large enough), we have
$ Vol(B_s)\geq C\cdot s^{n-\delta_1}$ for $s$ large enough. Similarly,
$\frac{1}{Vol(B_s)}\int_{B_s} R\geq \delta_2>0 $ implies $ Vol(B_s)\leq C\cdot s^{n-\delta_2}$.
}%%%%%%%%%%%%%%%%%%%%%%%%%%%%%%%%%%%%%%%%%%%%%%%%%%%%%%%%%%%%%%%%%%%%%%%
\\

A similar result to the case $a=0$ in Theorem 4 was proved by Cao and Zhou in \cite{CZ} (the case of volume lower bound).
The last statement concerning the sharp upper volume bound was proved before by S. Zhang in \cite{SZ}.
Moreover, Cao and Zhou \cite{CZ} and O. Munteanu \cite{Mu} proved that the upper bound $ Vol(B_s)\leq C\cdot s^n$ always holds for
{\it all} shrinking gradient Ricci solitons.

{\Rmk After the first version of this article has been posed on the website ArXiv.org,
B. Chow, P. Lu and B. Yang \cite{CLY}, based on the ideas and techniques of certain former works,
concluded a criterion for a shrinking soliton to have positive asymptotic volume ratio.
We remind the reader that, for non-flat gradient shrinking solitons which have nonnegative Ricci curvature,
Carrillo and Ni \cite{CN} have proved that they must have zero asymptotic volume ratio.
It generalizes a well-known result of G. Perelman \cite{Pe}.
}\\

From now on, we consider geodesic balls $B_r(x)$ whose center $x$ varies on $M$.
In general, or even for manifolds with fast decay curvature,
it is impossible to obtain a uniform volume bound of all $B_r(x)$ only from the
information of the aysmptotic volume ratio.
However, the following theorem shows that such bound does exist
if the manifold is an expanding soliton.

{\Thm%%%%%%%%%%%%%%%%%%%%%%%%%%%%%%%%%%%%%%%%%%%%%%%%%%%%%%%%%%%%%%%%%%%%%%%%%%%%%%%%%%%%%%%%%%%%%%
If a complete non-compact expanding gradient Ricci soliton $(M,g,f)$ satisfies $\lim_{s\rightarrow\infty} s^2\cdot |Sect|=0$,
then we have $$ Vol(B_{r}(x))\geq Cr^n$$ for all $r>0$ and $x\in M.$
We also have $ Vol(B_{r}(x))\leq Cr^n$ for all $r\leq \frac{s}{2}$ and for all $x\neq O.$
Moreover, if its asymptotic volume ratio exists,
then we have $$C^{-1}r^n \leq Vol(B_{r}(x)) \leq Cr^n$$ for all $r>0$ and $x\in M.$
}%%%%%%%%%%%%%%%%%%%%%%%%%%%%%%%%%%%%%%%%%%%%%%%%%%%%%%%%%%%%%%%%%%%%%%%%%%%%%%%%%%%%%%%%%%%%%%%%%

The following calculation is the crucial ingredient to achieve this theorem.
Given an expanding soliton $(M,g,f)$ with $\lim_{s\rightarrow\infty} s^2\cdot |Sect|=0$,
we consider a level set $\Sigma_a:=\{f=a\}$ of $f$.
The unit outer normal vector of $\Sigma_a$ shall be $\frac{-\nabla f}{|\nabla f|}$ and
the second fundamental form of $\Sigma_a$ is
$$II_{ij} =\left<\nabla_{e_i} \frac{-\nabla f}{|\nabla f|}, e_j \right>\\
=\left< -\frac{\nabla_{e_i}\nabla f}{|\nabla f|} + \frac{(e_i|\nabla f|)\nabla f}{|\nabla f|^2}, e_j\right>
= \frac{Hess (-f)_{ij}}{|\nabla f|},$$
for all $ i,j=1,\cdots,n-1$. Hence by Gauss equation, one obtain a curvature estimate of $\Sigma_a$.
We will use this to control the volume of a portion of $\Sigma_a$ in the following proof of Theorem 5.
We learned how to prove this theorem from Gilles Carron,
who has known it in mind and kindly shared it with us.

\begin{proof}[Proof of Theorem 5]
Step 1. We prove first that the lower bound estimate holds for all $x\neq O$ and $r=\frac{s}{2}:=\frac{1}{2}dist(O,x)$.
It suffices to show that, for $s$ large enough,
$B_{\frac{s}{2}}(x)$ contains a "cube" whose volume is at least $\bar{\delta} s^n$ for some $\bar{\delta}$ independent of $x$.
Let $f(x)=a$, $\Sigma_a:=\{f=a\}$ and $II=Hess(-f)/|\nabla f|$ be the second fundamental form of $\Sigma_a$.
Since $\|Hess(-f)-g\| \leq |Ric| \in o(s^{-2})$ implies $\|II-\frac{g}{|\nabla f|}\|\in o(s^{-3})$,
by Gauss equation and the fast decay of the curvature of $M$, we have $|Sect^\Sigma-\frac{1}{s^2}| \in o(s^{-2})$.
Hence for $s$ large enough, there exists an intrinsic ball $B_{\delta s}^\Sigma(x)\subset \Sigma_a$
such that $Vol^\Sigma (B_{\delta s}^\Sigma)\geq Cs^{n-1}$.

Furthermore, by using the one parameter family of diffeomorphisms
$\{\varphi_{t}: \Sigma_a\rightarrow \Sigma_{a+t}\}_{t\in (-\frac{s}{10},\frac{s}{10})}$,
which is generated by $\frac{\nabla f}{|\nabla f|^2}$,
we know the cube $Cube:=\{ y\in M | y=\varphi_t(B_{\delta s}^\Sigma(x)),t\in (-\frac{s}{10},\frac{s}{10})\}$
is contained in $B_\frac{s}{2}(x)$ and $Vol(Cube)\geq \bar{\delta} s^n$
for some constant $\bar{\delta}$ which is independent of $x$ whenever $s$ is large enough.

Step 2. For $r\leq\frac{s}{2}$, by using the Bishop-Gromov's comparison, we have
\begin{eqnarray*}
Vol(B_r(x))
&\geq& \left(Vol_{\frac{-C}{s^2}}\left(B_r\right)\slash Vol_{\frac{-C}{s^2}}\left(B_{\frac{s}{2}}\right)\right)\cdot Vol(B_{\frac{s}{2}}(x))\\
&\geq& \left(Vol_{\mathbb{R}^n}\left(B_r\right)\slash Vol_{\frac{-C}{s^2}}\left(B_{\frac{s}{2}}\right)\right)\cdot \bar{\delta}s^n\\
&\geq& Cr^n.
\end{eqnarray*}
The last inequality comes from
$$ Vol_{\frac{-C}{s^2}}(B_s)=C\int_0^s \left( {\frac{s}{\sqrt{C}}}  \sinh\left(\frac{\sqrt{C}}{s}t\right)\right)^{n-1} dt
\leq Cs^n,$$
where $Vol_{\frac{-C}{s^2}}$ is the volume functional of the hyperbolic space with $Ric=\frac{-C}{s^2}.$

Step 3. For every ball $B_r(x)$ with $r>\frac{s}{2}$ and $x\in M$, it must contain a ball $B_{\frac{r+s}{4}}(y)$
for some $y$ satisfying that $dist(O,y)=\frac{r+s}{2}$. By Step 1, we know that
$Vol(B_{\frac{r+s}{4}}(y))\geq C (\frac{r+s}{4})^n\geq C(\frac{r}{4})^n$. Hence $Vol(B_r(x))\geq Cr^n$.

Step 4. The upper bound can be derived by by Bishop-Gromov's comparison.
Since $Ric\geq -C\cdot s^{-2}$ on $B_r(x)$ with $r\leq \frac{s}{2}$, as in the last inequality of Step 2,
we have
   $$ Vol(B_r(x))
   \leq Vol_{\frac{-C}{r^2}}(B_r)
   =    C\int_0^r \left( {\frac{r}{\sqrt{C}}}  \sinh\left(\frac{\sqrt{C}}{r}t\right)\right)^{n-1} dt
   \leq Cr^n.$$

Therefore we have proved the first statement of the theorem.
Now suppose that we can control the upper bound of the volume ratio at infinity,
i.e., there exist two constants $C$ and $A$ such that $Vol(B_r(O))\leq Cr^n$ for all $s\geq A$.
It is easy to see that, for all $r>\frac{s}{2}$, $B_r(x)$ is contained in $B_{3r}(O)$ and
hence has an upper bound on its volume ratio.
\end{proof}

It was known that a Ricci-nonnegative expanding soliton must be diffeomorphic to $\mathbb{R}^n$ because
$-f$ is proper and strictly convex.
On the other hand, by using the same argument of W. Wylie in \cite{W} and our method developed in section 1,
one can prove that an expanding soliton with $Ric\geq -Cs^{-\varepsilon}$ must have finite fundamental group.
From Step 1 in the above proof of Theorem 5, we have the following topological information of the ends
of an expanding soliton.

{\Thm %%%%%%%%%%%%%%%%%%%%%%%%%%%%%%%%%%%%%%%%%%%%%%%%%%%%%%%%%%%%%%%%%%%%%%%%%%%%%%%%%%%%
Let $(M,g,f)$ be a complete non-compact expanding gradient soliton with $\lim_{s\rightarrow\infty} s^2\cdot |Sect|=0$.
Then each end of $M$ is diffeomorphic to $\mathbb{R}\times N^{n-1},$ where $N=\mathbb{S}^{n-1}/\Gamma$
is a metric quotient of the spherical space form.
}\\%%%%%%%%%%%%%%%%%%%%%%%%%%%%%%%%%%%%%%%%%%%%%%%%%%%%%%%%%%%%%%%%%%%%%%%%%%%%%%%%%%%%%%%%%%%%%%%%%%%%%

%%%%%%%%%%%%%%%%%%%%%%%%%%%%%%%%%%%%%%%%%%%%%%%%%%%%%%%%%%%%%%%%%%%%%%%%%%%%%%%%%%%%%%%%%%%%%%%%%%%%%%%%%%%%
\section{Tangent cones at inifinity of expanding gradient Ricci solitons}%%%%%%%%%%%%%%%%%%%%%%
%%%%%%%%%%%%%%%%%%%%%%%%%%%%%%%%%%%%%%%%%%%%%%%%%%%%%%%%%%%%%%%%%%%%%%%%%%%%%%%%%%%%%%%%%%%%%%%%%%%%%%%%%%%%

A tangent cone at infinity is a Cheeger-Gromov limit of a sequence of blow-down metrics with a fixing marked point.
Since we have a uniform estimate of volume lower bound from Theroem 5,
we can derive a lower bound of injectivity radius from the controlled sectional curvature (see \cite{CLY81} and \cite{CGT}.)
In this section, we prove that every tangent cone at infinity of $M$ is the Euclidean space $\mathbb{R}^n$
under some admissible conditions.

{\Thm
Let $(M,g,f)$ be a complete non-compact expanding gradient Ricci soliton which satisfies
$R_{ij}+\nabla_i\nabla_jf=- g_{ij}$ and $lim_{s\rightarrow\infty} s^2\cdot |Sect|=0$.
If $M$ is simply connected at infinity, has only one end and has dimension $n\geq 3$,
then every tangent cone at infinity of $M$ is the Euclidean space $\mathbb{R}^n$.
}

\begin{proof} %[Proof of Theorem 7]
Consider a tangent cone at infinity $M^\infty$, which is a Gromov-Hausdorff limit of
a sequence $(M,O,\widetilde{g}_k):=(M,O,\frac{1}{\lambda_k^2}g)$ with vertex $O$, where
$\lambda_k\rightarrow\infty$ as $k\rightarrow\infty$. Here we use a tilde to
emphasize that the metric is rescaled. Any arbitrary point $q\in M^{\infty },q\neq O$
and $dist_{\infty }(O,q)=r_{0},$ is associated with a sequence $q_{k}\rightarrow q,$
where $dist_{k}(O,q_{k})=\lambda_{k}r_{0}\rightarrow\infty$ as $k$ $\rightarrow\infty.$
By using our volume estimate in the previous section, Hamilton's compactness theorem
and Shi's estimate, the convergence is in fact in $C_{loc}^{\infty}$-topology.

Noting that $\left|\widetilde{\nabla}_i\widetilde{\nabla}_j f_k\right|=\left|(\widetilde{g}_k)_{ij}+\frac{1}{\lambda_k^2}(\widetilde{Ric}_k)_{ij}\right|$,
together with the estimates of the growth of $f$ and $\nabla f$ which are stated in Section 2,
we know that $f_k:=\frac{-f}{\lambda_k^2}$ converges in $C_{loc}^{\infty}$-topology to a function $f^\infty$
with $|\nabla f|=r$ on $M^\infty\setminus\{O\}$.
Moreover, $\nabla^\infty\nabla^\infty f^\infty=g^\infty$ and
$f^\infty(q)=\lim_{k\rightarrow\infty}\frac{-f}{\lambda_k^2}(q_k)=\frac{1}{2}r_0^2.$
Since $q$ was chosen arbitrarily, we have
$$f^\infty(x)=\frac{1}{2}r^2 \mbox{ \ \ and \ \ } g^\infty=Hess\left(\frac{1}{2}r^2\right)$$
where $r(x):=dist_\infty (O,x)$ and $x\in M^\infty\setminus\{O\}$.

In \cite{CC}, J. Cheeger and T. H. Colding have proven that $M^\infty\setminus\{O\}$
with $g^\infty=Hess(\frac{r^2}{2})$ must be a warped product manifold and $g^\infty=dr^2+kr^2\bar{g}$ for some $k>0$,
where $\bar{g}$ is the metric of $N:=\{x\in M^\infty | r(x)=1\}$.
In order to prove that $M^\infty$ is isometric to $\mathbb{R}^n$,
we only need to show that $N$ is the standard sphere with sectional curvature $k$.
(Because the standard metric on $\mathbb{R}^n$ can be written as
$g_{Eucl}=dr^2+Cr^2g_{\mathbb{S}^{n-1}(C)}$ for any given $C>0$ and
$g_{\mathbb{S}^{n-1}(C)}$ denotes the standard metric on sphere with constant sectional curvature $C$.)

Since $|\nabla r|\neq 0$, we can extend the normal coordinate $\{x^i\}_{i=2,\dots,n}$ around $p\in N$
to be a local coordinate $\{r,x^i\}_{i=2,\dots,n}$ in $M$ such that
$$(g_{ij})=\left(     \begin{array}{cccc}
                       1 & 0       & \cdots  &   0    \\
                       0 & g_{22}  & \cdots  &   g_{2n}    \\
                  \vdots & \vdots  & \ddots  & \vdots \\
                       0 & g_{n2}    & \cdots  & g_{nn} \\
                \end{array} \right)
     =\left(    \begin{array}{cccc}
                       1 &     0        & \cdots  &   0       \\
                       0 &   kr^2\bar{g}_{22}       & \cdots  &   kr^2\bar{g}_{2n}       \\
                  \vdots &   \vdots     & \ddots  & \vdots    \\
                       0 &   kr^2\bar{g}_{n2}        & \cdots  &  kr^2\bar{g}_{nn}     \\
                \end{array} \right).$$
Hence, for all $i,j=2,\dots,n$ and $i\neq j$, we have $\Gamma_{jj}^r(p)=-k$ and $\Gamma_{ij}^r(p)=0$.
Moreover, $\frac{\partial}{\partial x^j}( g(\frac{\partial}{\partial r},\frac{\partial}{\partial x^j}))=0$ implies that
$\Gamma_{jr}^j(p)=-\frac{1}{k}\Gamma_{jj}^r(p)=1$.
When $n\geq 3,$ we can compute the curvature of $N$ at $p$ by using
$$0=R_{ijj}^i=\bar{R}_{ijj}^i+\Gamma_{ir}^i\Gamma_{jj}^r=\bar{R}_{ijj}^i-k.$$
By the assumption that $M$ is simply connected at infinity, we know that
$N$ must be the standard sphere with all its sectional curvatures equal $k$.

\end{proof}

For the shrinking case, B. Chow, P. Lu and B. Yang \cite{CLY2} (by using a result of L. Ni and B. Wilking)
proved that a non-compact non-flat shrinking gradient Ricci soliton has at most quadratic scalar
curvature decay. Hence our theorem trivially holds for the shrinking case.
On the other hand, there exists a two-dimensional counter-example for the expanding case,
i.e. an expanding soliton which has faster-than-quadratic-decay curvature
and a tangent cone at infinity which is not an Euclidean plane.
Such soliton was constructed in \cite{CLN} by smoothly extending a cone manifold which had been
conceived in \cite{GHMS}.

%%%%%%%%%%%%%%%%%%%%%%%%%%%%%%%%%%%%%%%%%%%%%%%%%%%%%%%%%%%%%%

\end{document}